\documentclass[11pt,a4paper]{article}
\usepackage[dvips]{epsfig}
\usepackage{graphicx,latexsym, amsmath,amsfonts}
\usepackage{psfrag}
\usepackage{parskip}
\usepackage{graphicx}
\usepackage[]{epsfig}
\usepackage{amsmath}
\usepackage{amssymb}
\usepackage{verbatim}
\usepackage{url}
\usepackage{latexsym}
\usepackage{color}

\newtheorem{theorem}{Theorem}[section]

\begin{document}
\date{}
\title{\bf Effect of salinity and fish predation on zooplankton dynamics in Hooghly-Matla estuarine system, India}
\author{Ujjwal Roy$^{(1)}$,~S. Sarwardi$^{(2)}$\footnote{Author to whom all correspondence should be addressed},~N. C. Majee$^{(1)}$,~Santanu Ray$^{(3)}$\\
$^{(1)}$ Department of Mathematics, Jadavpur University\\ Kolkata-700 032, West Bengal, India\\
email:~\textcolor{blue}{ujjwal2ju@gmail.com}\\
$^{(2)}$ Department of Mathematics, Aliah University, IIA/27, New Town\\
Kolkata - 700 156, West Bengal, India\\
email:~\textcolor{blue}{s.sarwardi.math@aliah.ac.in}\\
$^{(3)}$Ecological Modelling Laboratory, Department of Zoology,
 Visva-\\Bharati University, Santiniketan-731 235, West Bengal, India.\\
email:~\textcolor{blue}{sray@visva-bharati.ac.in}}
\maketitle
\thispagestyle{empty}
 \begin{abstract}
The Hooghly-Matla estuarine complex is the unique estuarine system of the world. Nutrient from the litterfall enrich the adjacent estuary through tidal influence which in turn regulate the phytoplankton, zooplankton and fish population dynamics. Environmental factors regulate the biotic components of the system, among which salinity plays a leading role in the regulation of phytoplankton, zooplankton and fish dynamics of the estuary. In this article, a $PZF$ model is considered with Holling type-II response function. The present model considers  salinity based equations on plankton dynamics of the estuary. The interior equilibrium is considered as the most important equilibrium state of this model. The model equations are solved both analytically and numerically using the real data base of Hooghly-Matla estuarine system. The essential mathematical features of the present model have been analyzed thorough local and global stability and the bifurcations arising in some selected situations. A combination of set of values of the salinity of the estuary are identified that helped to determine the sustenance of fish population in the system. The ranges of salinity under which the system undergoes Hopf bifurcation are determined. Numerical illustrations are performed in order to validate the applicability of the model under consideration.  \end{abstract}


\section{Introduction}
The delta of Hooghly-Matla estuarine system is networked by seven major river  along with their tributaries and creeks. This deltaic system harbours luxuriant mangroves and constitutes Sundarban mangrove ecosystem. The mangroves are the major resources of detritus and nutrients to the adjacent estuary that make up favourable habitat for the growth of shell fish and fin fish (cf. Mandal et al., (2009)). This estuary act as route and refuge areas for a variety of migratory fish species. The estuary supports 53 species of pelagic fish belonging to 27 families and 124 species under 49 families of demersal fish (cf. Hussain and Acharya, (1994)). Fishery in this estuarine water contributes a part  in the economy of the state of West Bengal, India. The organic matter that passes from litterfall to the adjacent estuary supports both the grazing and detritus food chains of this lotic system.

Among the chemical components studied so far, the salinity plays a crucial role in the abundance (cf. Bhunia, (1979)) and dynamics of zooplankton of the estuary (cf. Ghosh, (2001), Ketchum, (1951)). Because, this community lies in the middle of grazing food chain; phytoplankton-zooplankton-fish ($PZF$). In addition, perturbation  to this trophic level may trigger imbalance in the food chain which in turn affect the phytoplankton and fish community. The dynamics of salinity is season dependent. During  monsoon and early post monsoon, huge fresh water enters in the estuary from the upstream resulting in the lowering of salinity. In the pre-monsoon, fresh water runoff from the upstream becomes very less and due to tidal influence of the adjacent  Bay of Bengal, the salinity increases. Throughout the year, a gradient of  salinity is observed between upstream and downstream area of the estuary (cf. Mandal et al., (2009)).

The abundance of different species of zooplankton varies according to the salinity of the estuary throughout the year, as a result the grazing rate also changes with seasons.  In $PZF$ system, the grazing rate of zooplankton is one of the most sensitive parameter in Hooghly-Matla estuarine system  as the dynamics of zooplankton also depends on the lower trophic level of phytoplankton as well as fish population (cf. Mandal et al., (2009)). Moreover, fish predation exhibits top-down effect on zooplankton community. Therefore, in the $PZF$ system, two important parameters that shape up the zooplankton dynamics of the estuary are salinity dependent grazing rate of zooplankton and fish predation rate (cf. Dube et al., (2010)).

There have been only a few models based on the effects of salinity on plankton dynamics. Few studies have been done in the relationship between salinity levels to the types of species that can occur in algal blooms (cf. Griffin et al., (2001), Marcarelli et al., (2006), Quinlan and Phlips, (2007)). Many $PZF$ model are constructed on estuarine system (cf. Ray et al., (2001); Cottingham, (2004); Dube and Jayaraman, (2008); Dube et al., (2010)) and top-down effect (cf. Morozov et al., (2005); Irigoien et al., (2005); Calbet and Saiz, (2007)) but none of the models have considered the salinity dependent grazing rate of zooplankton which plays an important rule in the zooplankton dynamics of estuary. The present account deals with a $PZF$ model, where salinity dependent grazing rate of zooplankton is taken into consideration along with the variation of upstream and downstream salinity.

\subsection{Mathematical model formulation}
Let $P$, $Z$, $F$ denote the populations of phytoplankton, zooplankton and fish respectively.
In this present $PZF$ model, light and temperature dependant photosynthesis rate of phytoplankton and salinity induced grazing rate of zooplankton have been incorporated into the model proposed by Mandal et al., (2011). The modified model under consideration is as follows:
\begin{eqnarray}\left \{ \begin{array}{lll}
\frac{dP}{dt} = m_{1}P\Bigl(1-\frac{P}{k_{P}}\Bigr)-g_{s}\frac{PZ}{P+k_{Z}}\\
\frac{dZ}{dt} = a g_{s}\frac{PZ}{P+k_{Z}}-g_{f}\frac{ZF}{Z+k_{F}}-m_{2}Z\\
\frac{dF}{dt} = g_{f}\frac{ZF}{Z+k_{F}}-m_{3}F.\end{array}\right.
 \label{eq1}\end{eqnarray}


\text{The salinity induced grazing of zooplankton is}~~$g_s=\delta g_Z,~\text{where}~~  \delta =\frac{s_{u}}{s_{u}-s_{d}},$ the dilution factor and $m_{1},$ ~$s_{d},$~$s_{u},$~$g_{Z},$~$k_{Z},$~$k_{P}$ are
net growth rate of phytoplankton, downstream  salinity, upstream salinity, grazing rate of zooplankton, half saturation constant of zooplankton on phytoplankton grazing by zooplankton and half saturation constant grazing by phytoplankton respectively. Since estuary is a transition zone of river and sea, so there is always fluctuation of salinity throughout the year, which is due to dilution by upstream river water and/or mixing by downstream sea tidal water. $\delta$  is calculated by following the equation of Ketchum, (1951). Besides grazing the abundance of zooplankton is also dependent on loss due to $E_{zoo}=Z E_{zo} ,$ respiration $R_{zoo}=Z r_{zo},$ fish predation $F_p=Z r_{rf}$, and mortality $M_{zoo}=Z M_{Z}$, $R_{zoo}$ is governed by $r_{zo}$.

\text{The net mortality rate of zooplankton is }~$m_2=(E_{zo}+r_{zo} + r_{fp} +M_{Z}).$ The abundance of $F$ is governed by many processes in the estuary. Fish predation $F_{p}$ on $Z$ follows Michaelis-Menten kinetics (Holling type II) which enriches the fish pool of the estuarine system.
$F_{p}$ depends on $k_{Z}$ and $g_f$. $F$ population is reduced by mortality rate of fish $M_{f}$, respiration $Fr_f$, harvest by fishing $FH_{cf}$ and excretion of fish $FE_f$. $FH_{f}$ are controlled by $R_f$ and $H_f$ respectively. Respiration rate of  Fish and $H_{f}$ Harvest rate of Carnivorous Fish respectively.

\text{The net mortality rate of fish is}~~$m_3=(E_{f} +M_{f}+r_{f}+H_{f}),$
where $r_{fp},$~$r_{zo}$~and ~$M_{rz}$ are fish predation rate, zooplankton respiratory rate and mortality rate respectively.
\subsection{Existence and positive invariance}
 Letting $ X \equiv(P,Z,F)^T$,$ ~f: \mathbb{R}^3\rightarrow \mathbb{R}^3, f=(f_{1},f_{2}, f_{3})^T$, the system (\ref{eq1}) can be rewritten as $\dot{X}=f(X)$. Here $f_{i}\in C^\infty (\mathbb{R})$ for $i=1,2,3,$ where $f_{1}= m_{1}P(1-\frac{P}{k_{P}})-g_{s}\frac{PZ}{P+k_{Z}}$, $f_{2}= {a} g_{s}\frac{PZ}{P+k_{Z}}-m_{2}Z-\frac{ZF}{Z+k_{F}}$, $f_{3}= g_{f}\frac{ZF}{Z+k_{F}}-m_{3}F$. Since the vector function $f$ is a smooth function of the variables $P$, $Z$ and $F$ in the positive octant $ \Omega^0 =\{(P,Z,F);P>0,Z>0,F>0\},$ the local existence and uniqueness of the system hold.
\subsection{Boundedness of the system}
Boundedness of a system guarantees its biological validity. The following theorem establishes the uniform boundedness of the system (\ref{eq1}).

\textbf{Theorem 1.}  \emph{All the solutions of the system (\ref{eq1}) which start in $\mathbb{R}_{+}^3$ are uniformly bounded.}

\textbf{Proof.} Let $(P(t), Z(t), F(t))$ be any solution of the system with positive initial conditions.
From the real ecological field study one can consider ~ $\max \{ ~ m_{2},~ m_3\} < 1.$

Now let us we define the function $X =aP+Z+F.$ ~The time derivative of X gives
\begin{eqnarray}
\frac{d X}{d t} &=&a\frac{d P}{d t} + \frac{d Z}{d x} + \frac{d F}{d t},\nonumber\\
&=& am_{1}P(1-\frac{P}{k_{P}})-ag_{s}\frac{P Z}{P+k_{Z}} +a g_{s}\frac{P Z}{P+k_{Z}}-g_{f}\frac{Z F}{Z+k_{F}}-m_{2}Z + g_{f}\frac{Z F}{Z+k_{F}}-m_{3}F \nonumber\\
&=&am_{1}P(1-\frac{P}{k_{P}}) - m_2 Z - m_3 F.\end{eqnarray}
Now for each $v > 0$, we have
\begin{eqnarray}
\frac{d X}{d t} + v X &=&am_{1}P(1-\frac{P}{k_{P}}) - m_2 Z - m_3 F + v(aP+Z+F)\nonumber\\
&=& am_{1}P(1-\frac{P}{k_{P}})+ a v P +(v-m_2)Z+(v-m_3)F\nonumber\\
&\leq&am_{1}(1-\frac{P}{k_{P}})+ a v P,~ \text{if}~ v \leq  ~\min \{ ~ m_{2},~ m_3\},\nonumber\\
~&\leq &\frac{a (m_{1}+v)^2k_{P}}{4}.
\end{eqnarray}
Using the theory of differential inequality (cf. Birkoff and Rota, (1982); Sarwardi et al. (2011)), one can easily obtain
\begin{eqnarray}
\limsup_{t\rightarrow +\infty}X(t)\leq \frac{a (m_{1}+v)^2k_{P}}{4}=\rho.
\end{eqnarray}
Therefore, all the solutions of the system (\ref{eq1}) enter into the compact region $\mathcal{B}=\Bigl\{(P,Z,F)\in \mathbb{R}_{+}^3 : aP+Z+F \leq \frac{a (m_{1}+v)^2k_{P}}{4}\Bigr\},$
which completes the proof.

\subsection{Equilibria and their feasibility}
To determine the equilibrium points of the system of(\ref{eq1}) we put $\dot{P}=\dot{Z}=\dot{F}=0$.
The equilibria of the system(\ref{eq1}) are
(i) the null equilibrium $E_{0}=(0,0,0);$\\
(ii) the axial equilibrium $E_{1}=(k_{P},0,0);$\\
(iii) the boundary equilibrium $E_{2}=(\frac{m_{2}k_{Z}}{(a g_{s} -m_{2})},~\frac{(a m_{1} k_{Z})(ak_{P}g_{s}-k_{P}m_{2}-m_{2}k_{Z})}{(a g_{s}-m_{2})^2 k_{P}},0);$~and\\
(iv) the interior equilibria ~~$E_{*}^i=(P_i^*, Z^*, F_i^*);~i = 1,2,$

where $P_i^*$=$\frac{-m_1(g_f-m_3)(k_{P}-k_Z)\pm\sqrt{(m_1(g_f-m_3)(k_{P}-k_Z))^2-4{m_1(m_3-g_f)}(k_{P}k_Zm_1(g_f-m_3)-m_3k_fg_sk_P) }}{2m_1(m_3-g_f)},$ \\

$Z_i^*$=$\frac{m_3k_F}{(g_f-m_3)}$ and $F_i^* = \frac{m_1k_Fa(P^*-k_P)}{m_3k_P}+\frac{k_F(ag_s-m_2)}{(g_f-m_3)}.$

\subsection{Existence of planer and interior equilibria and their stability}
Equilibria of the model (\ref{eq1}) can be obtained by solving
$\dot{P}=\dot{Z}=\dot{F}=0$. Though the system (\ref{eq1}) has several
non-negative steady state but we have mentioned the existence of interior equilibrium point only for its biological importance.\\
The interior equilibrium point $E_{*}=(P^*,Z^*,F^*)$ of the system (\ref{eq1}) exists if $ (g_f-m_2)>0$ . When these conditions are satisfied, we have $P^*$ as the positive root of the following quadratic equation:
\begin{equation}
A_0 {P^*}^2 - A_1 P^* - A_2=0,\label{eq3}
\end{equation}
where
\begin{eqnarray}
 A_0 =-m_1(g_f-m_3),~ \, A_1 =(g_fm_1-m_3m_1)(k_{P}-k_Z), ~\,
A_2 = k_{P}k_Zm_1(g_f-m_3)-m_3k_fg_sk.\nonumber
\end{eqnarray}
Since $g_f>m_3,$ $A_0$ is positive. Therefore, one positive root of (\ref{eq3}) can be found as
\begin{equation}
P^*=\frac{1}{2A_0}[A_1+\sqrt{A_1^2+4A_0A_2}],
\end{equation}

which exists if $A_2>0$ and $A_0>0$.
Therefore, the equation (\ref{eq3}) has a unique
positive solution if
$g_{f} > m_{3}.$
Once we get the unique positive solution of $P^*$ from equation (\ref{eq3}), it is easy to obtain the other components of the interior equilibrium $E_*.$

\textbf{Feasibility:}
It is clear that the axial equilibrium is feasible. The interior equilibrium is feasible if the condition (i) $ (g_{f} - m_{3}) > 0 $ and (ii) $(k_{P}-k_z) < 0 $ hold good.
\subsection{Local stability of equilibrium}
In order to find out about the stability of the equilibrium points we need to linearize the system (\ref{eq1}). Since the axial and boundary equilibria are less important on this system (\ref{eq1}). That's why we don't go for detailed analysis of the system around the equilibria other than interior equilibrium. The stability of interior equilibrium, we find the Jacobian matrix. For this purpose we define the system (\ref{eq1}) as follows:
\begin{equation}
\frac{dX}{dt}= F(X),
\end{equation}
where
$f(X) = \left[
\begin{array}{c}
m_{1}P(1-\frac{P}{k_{P}})-g_{s}\frac{PZ}{P+k_{Z}}\\
   a g_{s}\frac{PZ}{P+k_{Z}}-g_{f}\frac{ZF}{Z+k_{F}}-m_{2}Z \\
   g_{f}\frac{ZF}{Z+k_{F}}-m_{3}F
  \end{array}
\right],$ ~ \text{and} ~ $Y=\left[\begin{array}{c}
P\\
Z\\
F
\end{array}
\right].$

The variational matrix of the system at any arbitrary point ( $\tilde{P}$,$\tilde{Z}$,$\tilde{F}$ ) is

$\tilde{J}=\frac{\partial f}{\partial X} \mid_{(\tilde{P},\tilde{Z},\tilde{F})} =\left[
\begin{array}{ccc} \displaystyle \frac{2 g_s \tilde{Z}}{\tilde{P}+k_Z}-m_1-\frac{g_s \tilde{Z} k_{Z}}{(\tilde{P} +k_{Z})^2}
&\displaystyle -\frac{g_s \tilde{P} }{(\tilde{P} + k_{Z})} &\displaystyle 0 \\
 \displaystyle \frac{ag_s k_{Z} \tilde{Z}}{(\tilde{P} +k_{Z})^2} &\displaystyle \frac{g_f \tilde{F} }{(\tilde{Z} + k_F)}-\frac{g_{f} \tilde{F} k_{F}}{(\tilde{Z} + k_{F})^2} & \displaystyle -\frac{g_{f} \tilde{Z}}{(\tilde{Z} + k_{F})} \\
 \displaystyle 0 & \displaystyle \frac{g_{f} \tilde{F} k_{F}}{(\tilde{Z} + k_{F})^2} &\displaystyle 0 \\
  \end{array}
\right].$

Therefore, the Jacobian matrix for the system (\ref{eq1}) at $E_*$ is given by

$$J_{*} =\left[
\begin{array}{ccc} \displaystyle R_1-R_2-m_1
&\displaystyle -K_{1} &\displaystyle 0 \\
 \displaystyle R_2 &\displaystyle K_{2}-R_3 & \displaystyle -m_{3} \\
 \displaystyle 0 & \displaystyle R_{3} &\displaystyle 0
  \end{array}
\right],$$

where
$R_{1}=\frac{2 g_s Z^*}{(P^* +k_Z)},~
R_{2}=\frac{g_s Z^* k_Z}{(P^* +k_Z)^2},~
R_{3}=\frac{g_f F^* k_F}{(Z^* + k_F)^2},~
K_{1}=\frac{g_s P^*}{(P^* + k_Z)},~
K_{2} =\frac{g_f F^* }{(Z^* + k_F)}.$

The characteristic equation corresponding to the Jacobian matrix $J_*$ is
\begin{equation}
\Delta(\lambda)=\lambda^3 +  D_1\lambda^2 + D_2 \lambda + D_3 = 0.
\end{equation}

For the stability of the solutions of the system (\ref{eq1}) all of the roots of the characteristic equation for the Jacobian matrix at $E_{*}$ ie. $J_{*}$ should have negative real parts. This can be achieved without
actually solving for all the roots of the characteristic equation, but by applying the Routh-Hurwitz condition for negative real parts of the characteristic roots. According to the Routh-Hurwitz conditions, the solutions of (\ref{eq1}) will be asymptotically stable if the conditions $ D_1 > 0,~~  D_1D_2 >D_3,~  D_3 > 0$ are satisfied.

\section{Global Stability}

%

Now we examine the global stability issue.
Here we describe the general method described by Li and Muldowney, (1996) to show an n-dimensional autonomous dynamical system
$f:D\rightarrow \mathbb{R}^n, D\subset \mathbb{R}^n$,
an open and simply connected set and  $f\in C^{1}(D)$, where the dynamical system is given by
\begin{equation}
\frac{dx}{dt}= F(x)\label{eq4}
\end{equation}
is globally stable under certain parametric conditions (cf. Jin and Haque, (2005);  Haque et al., (2009)  and Buonomo et al., (2008))\\
Now we assume the following incidents: \\
(A1) The autonomous dynamical system has a unique interior equilibrium point x$^*$ in D.\\
(A2) The domain D is simply connected.\\
(A3) There is a compact absorbing set $\Omega\subset$ D.\\
Now we give a definition due to Haque et al., (2005)

\textbf{Definition 1.} The unique equilibrium point x$^*$ of the dynamical system (\ref{eq4}) is globally asymptotically stable in the domain D if it is locally asymptotically stable and all the trajectories in D converges to its interior equilibrium point x$^*$.

Let J=(J$_{ij})_{3x3}$ be the variational matrix of the system (\ref{eq4}) and J$^{[2]}$ be second additive compound matrix with order $^{n}C_{2} \times$ $^{n}C_{2}.$
In particular for n=3 we can write \\

$J^{[2]} =\frac{{\partial f}^{[2]}}{\partial X}=\left[
\begin{array}{ccc} V_{11}+V_{22}
& V_{23} & -V_{13} \\
 V_{32} & V_{11}+V_{33} & V_{12}\\
  -V_{31} & V_{21} & V_{22}+V_{33}  \\
  \end{array}
\right]$
Let $M(x) \in C^{1}(D)$, be an $^{n}C_{2}$ $\times$ $^{n}C_{2}$ matrix valued function. Moreover, we consider B as a matrix such that B=M$_{f} M^{-1}+MJ^{[2]}M^{-1}$ where the matrix M$_{f}$ is represented by
\begin{equation}
 (M_{ij}(x))_{f}= {\biggr(\frac{\partial M_{ij}}{\partial x}\biggr)}^t.f(x)=\nabla M_{ij}.f(x).\label{eq5}\\
\end{equation}
Again due to Jr. Martin et al., (1974), we consider the Lozinski$\check{i}$ measure $\Gamma$ of $B$ with respect to a vector norm $|.|$ in $\mathbb{R}^n$ . $N=$ $^{n}C_{2}$, we have
\begin{equation}
\Gamma (B)=\lim_{h \rightarrow 0^+}\frac{|I+h B|-1}{h}. \notag{}
\end{equation}
 Li and Muldowney (1996) showed that the system (\ref{eq4}) will be globally stable if the conditions (A1), (A2) and (A3) together with
\begin{equation}
\limsup_{t \rightarrow \infty} \sup \frac{1}{t}\displaystyle\int_{0}^{t}\Gamma (B(x(s,x_{0})))ds<0\label{eq6}
\end{equation} hold simultaneously.

The above stated condition not only assures that there are no orbits (i.e., homoclinic orbits, heteroclinic cycles and periodic orbits) which gives rise to a simple closed rectifiable curve in D, invariant for the system (\ref{eq4}) but also it is a robust Bendixson criterion for (\ref{eq5}) in addition with (\ref{eq6}) is sufficient for local stability the system (\ref{eq4}) around $E_{*}.$
Now we use the above discussion to show that our system (\ref{eq1}) is globally stable around its interior equilibrium.
We make the following transformation of the variables
$x\rightarrow P^{-1}$,~$y\rightarrow Z$, ~$z\rightarrow F$, for which the autonomous system (\ref{eq1}) is transformed to the following one:
\begin{equation}
\frac{dX}{dt}= f(X),\label{eq7}
\end{equation}
where
$f(X) =\left[
\begin{array}{c}
 \frac{m_1}{x}(1-\frac{1}{xk_P})-\frac{g_s y}{1+xk_Z}-x^2\\
   \frac{ag_sy}{1+xk_Z}-\frac{g_f yz}{y+k_F}-m_2y \\
   \frac{g_f yz}{y+k_F}-m_3 z\\
  \end{array}
\right]$ ~ and ~ $X=\left[
\begin{array}{c}
x\\
y\\
z\\
\end{array}
\right].$
The variational matrix of the system (\ref{eq7}) can be written as

$V =\frac{\partial f}{\partial X}=\left[
\begin{array}{ccc} \displaystyle \frac{x g_sy(2+xk_Z)}{(1 +xk_Z)^2}-m_1
&\displaystyle \frac{g_s x^2 }{(1 + xk_Z)} &\displaystyle 0 \\
 \displaystyle -\frac{ag_s y k_Z}{(1 +xk_Z)^2} &\displaystyle \frac{ag_s}{(1 + xk_Z)}-\frac{a g_s y}{(1 + xk_Z)}+\frac{g_fyz}{(y+k_F)^2} & \displaystyle -\frac{g_fy}{(y + k_F)} \\
 \displaystyle 0 & \displaystyle \frac{g_f Z k_F}{(y + k_F)^2} &\displaystyle 0 \\
  \end{array}
\right].$

If $V^{[2]}$ be the second additive compound matrix of V then due to Buonomo et al. (2008) we can write

$V^{[2]} =\left[
\begin{array}{ccc} \displaystyle a_{11}
&\displaystyle -\frac{g_f y}{(y + k_F)}  &\displaystyle 0 \\
 \displaystyle \frac{g_f Z k_F}{(y + k_F)^2} &\displaystyle a_{22}
 & \displaystyle \frac{g_s x^2 }{(1 + xk_Z)} \\
\displaystyle 0 & \displaystyle -\frac{a g_s y k_Z}{(1 +xk_Z)^2} &\displaystyle a_{33}
  \end{array}
\right],$
where
\begin{eqnarray} a_{11}&=&\frac{x g_sy(2+xk_Z)}{(1 +xk_Z)^2}-m_1+\frac{g_f y z }{(y + k_F)^2},\nonumber\\
a_{22}&=&\frac{x g_sy(2+xk_Z)}{(1 +xk_Z)^2}-m_1,\nonumber\\
a_{33}&=&\frac{g_fyz}{(y+k_F)^2} .
\nonumber\end{eqnarray}
We consider $P(X)$ in $C^{1}(D)$ in a way that $P=$ \text{diag} $\bigl(1, \frac{y}{z}, \frac{y x^2}{z}\bigr)$
then we have $P^{-1} = \text{diag} \bigl(1, \frac{z}{y}, \frac{z}{yx^2}\bigr)$,

$P_f P^{-1} = \text{diag} \{0, \frac{\dot{y}}{y}- \frac{\dot{z}}{z}, ~ \frac{2\dot{x}}{x} + \frac{\dot{y}}{y}-\frac{\dot{z}}{z}\}$ \text{and}

$PV^{[2]}P^{-1} =\left[
\begin{array}{ccc} a_{11}
& -\frac{g_fy}{(y + k_F)}  & 0 \\
  \frac{g_f  k_F z}{(y + k_F)^2} & a_{22}
 & \frac{g_sx^2}{(1+ x k_Z)} \\
 0 &  - \frac{a g_s y k_Z}{(1 +xk_Z)^2} & a_{33}
  \end{array}
\right].$
After some algebraic calculation, we get

$B = P_f P^{-1} + P V^{[2]} P^{-1}=\left[\begin{array}{ccc} B_{11}&B_{12} \\
  B_{21} &B_{22}  \\
  \end{array}
\right]$,

where

 $B_{11}=a_{11}= \frac{x g_sy(2+xk_Z)}{(1 +xk_Z)^2}-m_1+\frac{g_f y z }{(y + k_F)^2}$ ,\\
 $B_{12}=
\left[
  \begin{array}{ccc} -\frac{g_f y}{y+k_F}  & 0  \\
  \end{array}
\right]$,
 $B_{21}=
\left[
  \begin{array}{ccc} \frac{g_f z k_F}{(y + k_F)^2}  &  0  \\
  \end{array}
\right]^T$,
 $B_{22}= \left[
\begin{array}{ccc}
  a_{22}+ \frac{\dot{y}}{y}-\frac{\dot{z}}{z} & \frac{g_sx^2}{(1 + x k_Z)} \\
-  \frac{a g_s y k_Z}{(1 + x k_Z)^2} & a_{33}+\frac{2\dot{x}}{x}+\frac{\dot{y}}{y}-\frac{\dot{z}}{z}  \\
  \end{array}
\right]$.

Now let us define the following vector norm in $R^3$ as
$|(u,v,w)|=$max $\big( |u|,|v|+|w|)$, where $|(u,v,w)|$ is the vector norm in $R^3$ and it is denoted by $\Gamma$, the Lo\v{k}inski measure with respect to this norm.
Therefore, $\Gamma(B)\leq$sup$\{l_1,l_2\}$, where $l_i=\Gamma_1($B$_{ii})+|$B$_{ij}|$ for i$=$1, 2 and i$\neq$ j, here $|$B$_{12}|$ and $|$B$_{21}|$ are the matrix norms with respect to the L$^1$ vector norm and $\Gamma_1$ is the Lo\v{k}inski measure with respect to that norm. Therefore, we can easily obtain the following terms:
$\Gamma_1(B_{11})=\frac{x g_s y(2 + x k_Z)}{(1 +x k_Z)^2}-m_1+\frac{g_f y z }{(y + k_F)^2},
|B_{12}|= \max \{ \frac{g_f y}{y+k_F}, \hspace{0.2cm} 0 \},\\
|B_{21}|= \max \{ \frac{g_{f}  z  k_{F}}{(y + k_F)^2},\hspace{0.2cm} 0\},\\
\Gamma_1(B_{22})= \frac{\dot{y}}{y}-\frac{\dot{z}}{z} + \max\left\{\frac{g_s y (2x+x^2 k_Z-a k_Z)}{(1+x k_Z)^2},\hspace{0.2cm} \frac{g_s x^2}{(1+ xk_Z)}+\frac{g_f y z}{(y+k_F)^2}+\frac{2\dot{x}}{x}) \right \}.$
Using the system equation (\ref{eq7}),
 we have
 $l_1 = \frac{x g_s y(2 + x k_Z)}{(1 +x k_Z)^2}-m_1+\frac{g_f y z }{(y + k_F)^2} +\frac{g_f y}{y+k_F},$~~
$l_2 = \frac{\dot{y}}{y}-\frac{\dot{z}}{z} -2m_{1}(1-\frac{1}{xk_P})+ \frac{g_s y x}{1+xk_Z}+\frac{g_f y k_F}{(y+ k_F)^2}.$\\

Now from the expression of $l_1$ and $l_2$ one can obtain

 $\Gamma(B)=\max\biggr\{ \frac{x g_s y(2 + x k_Z)}{(1 +x k_Z)^2}-m_1+\frac{g_f y z }{(y + k_F)^2} +\frac{g_f y}{y+k_F},  \frac{\dot{y}}{y}-\frac{\dot{z}}{z} -2m_{1}(1-\frac{1}{xk_P})+ \frac{g_s y x}{1+xk_Z}+\frac{g_f y k_F}{(y+ k_F)^2}\biggr\}.$\\
 $~~~~~~~~=\frac{\dot{y}}{y}-\frac{\dot{z}}{z} -2m_{1}(1-\frac{1}{xk_P})+ \frac{g_s y x}{1+xk_Z}+\frac{g_f y k_F}{(y+ k_F)^2},$ if
 $1<x k_P<2.$

Thus,

$ \Gamma(B)\leq\frac{\dot{y}}{y}-\biggr[\frac{2g_s}{\rho+k_F}-(\frac{g_s}{k_Z}+\frac{g_f}{k_F})\rho-(m_1+m_3)\biggr]\hspace{4.2cm}\\
\,~~~~ ~~~~\leq\frac{\dot{y}}{y} -\mu,$ where $\mu=\frac{2g_s}{\rho+k_F}-(\frac{g_s}{k_Z}+\frac{g_f}{k_F})\rho-(m_1+m_3)>0.$

Taking average value of $\Gamma(B)$ within the time scale $0$ and $t$ one can have
\begin{equation}
\hspace{0.2cm} \frac{1}{t}\int^t_0\Gamma(B)ds\leq\frac{1}{t}log\frac{y(t)}{y(0)}-\mu, \notag{}\\
\text{which gives}
\end{equation}
\begin{equation}
\lim_{t \rightarrow \infty}sup\hspace{0.2cm} sup \hspace{0.2cm}\frac{1}{t}\int^t_0\Gamma(B(x(s,x_0)))ds<-\mu<0, \label{eq8}
\end{equation}
whenever $1>\frac{1}{x k_{P}}>\frac{1}{2}.$

Now we are in a position to state the following theorem:
\begin{theorem}The model system (\ref{eq1}) is globally asymptotically stable around its interior equilibrium if $\mu>0$,
 i.e $\frac{2g_s}{\rho+k_F}-(\frac{g_s}{k_Z}+\frac{g_f}{k_F})\rho-(m_1+m_3)>0.$
 \end{theorem}

 \subsection{Numerical simulation}
For the purpose of making qualitative analysis, numerical simulations have been carried out by making use of MATHLAB-R2010a and Maple-18. In Table 1, we summarize the used parameters with their admissible values and their biological interpretations. All these results have also been verified by numerical simulations,of which we report some in the figures. In order to compare our model with the corresponding model. we have run simulations using the standard MATLAB differential equations integrator for the Runge-Kutta method, i.e, MATLAB routine ODE45. The numerical experiments performed on the system (\ref{eq1}) using the experimental data taken from different sources to confirm our theoretical findings. Specially we focus here on the following values of the system parameters as shown in the table-I. The problem described by the system (1) is well posed. The P, Z, F axes are invariant under the flow of the governing system. In the present system the total environmental population is bounded above (cf. Subsection 1.3). Therefore, any solution starting in the interior of the first octant never leaves it. This mathematical fact is consistent with the biologically well behaved system and is common in the modern research work on mathematical biology.
 We took a set of parameter values: $m_1=0.6,\, m_2=0.0698,\, m_3=0.324,\, a=0.8,\, g_s=0.75,\, g_f=0.6894, \,s_d=12.30,\, s_u=8.23,\, k_P=12,\, k_Z=38,\, k_F=10.1.$ For this set of parameter values the system possessed a unique equilibrium point  $E_* (1.809, 8.964, 3.112)$. We have established the sufficient conditions for global stability of the coexistence equilibrium. The coexistence coexistence equilibrium point $E_{*}$ has been found through numerical simulations whose global asymptotical stability has been depicted in Figures 1(a) and Figure 2(a). It is observed from the Figure 1:(b)-(c) and Figure 2(b)-(c) that the fish population gradually decreases as the upstream salinity increases from $s_u=8.33$ to $s_u = 8.51$ while the downstream salinity is keeping fixed at $s_d=12.30.$ For higher upstream salinity ($s_u=8.51$) Figure 1(c) in the estuary is detrimental to the present zooplankton population. Zooplankton in this particular region are stenohaline thus the fish population cannot persist in such high salinity and the population gradually decreased to zero. Therefore, the salinity has an important role on regulating fish population in the present estuarine system. Representative numerical simulations are shown in Figures: 3--5. In Figure 3(a), the solutions plots describing the Hopf bifurcation and Figure 3(b) exhibits the corresponding phase space diagram. The present dynamical system confirms
the existence of chaotic nature which has been presented in Figure 4:(a)-(b). Figure 5:(a)-(f) shows different phase space diagram for different values of upstream salinity within the range 4.50--7.50, while downstream salinity has been kept fixed at $s_d=10.30$ and the remaining parameter values are same as taken in Figure 1. It is observed the catastrophic change in the behaviour of the solution plots as the parameter $s_u$ varies from 4.50 to 8.0 (cf. Figure 5:(a)-(f)). More specifically, Figure 5(a) shows periodic orbit with period 4; Figure 5(b) shows the chaotic orbit; Figure 5(c) shows periodic orbit with period 6; Figure 5(d) shows the period doubling bifurcation; Figure 5(e) demonstrates a limit cycle and finally Figure 5(f) ensures the stability of the system around an interior equilibrium point $(2.622, 8.878, 10.12)$ of the system (\ref{eq1}).
\begin{table}
\centering
\caption{The set of parameter values and their biological terminologies}
\vskip 0.5cm
\begin{tabular}{ |c| c| c | c |c| c } \hline\hline
Descriptions & \parbox[t]{0.6in}{Fixed Parameters}& Values &
\parbox[t]{0.5in}{Units}& \parbox[t]{0.80in}{References} \\[3ex]
\hline\hline
%

\parbox[t]{1.6in}{Growth rate of phytoplankton} &${ m_{1} }$&
0.6 &$day^{-1}$&
\parbox[t]{1.55in}{Mandal et al., (2012)}\\[2.5ex] \hline

\parbox[t]{1.6in}{Half saturation constant for
nutrient uptake by phytoplankton} &${k_{P}}$ & \parbox[t]{0.3in}{$12$}&
\parbox[t]{1.1in}{dimensionless}&
\parbox[t]{1.55in}{Mandal et al., (2012)} \\[1.5ex] \hline

\parbox[t]{1.6in}{ Grazing rate of zooplankton}& ${g_{z}}$ & \parbox[t]{0.3in}{$0.75$} &$day^{-1}$ &\parbox[t] {1.6in}{Ray and Straskraba, (2001)} \\[1.5 ex]
\hline

 \parbox[t]{1.6in}{Optimal light intensity} & ${E_{zo}}$ & \parbox[t]{0.3in}{ $0.04$} & $day^{-1}$ &
\parbox[t]{1.55in}{Ray et al., (2001)}  \\[0.5ex]
\hline
 \parbox[t]{1.6in}{ Respiration rate of zooplankton }& ${ r_{zo} }$  &\parbox[t]{0.3in}{ $0.0153$}& {$day^{-1}$}
 &\parbox[t]{1.55in}{Ray et al., (2001)} \\[2.5ex] \hline

\parbox[t]{1.6in}{Fish predation} & ${ r_{fp} }$  &\parbox[t]{0.3in}{$0.2$}&  $day^{-1}$
 &\parbox[t]{1.55in}{Ray et al., (2001)} \\[2.5ex]
 \hline
 \parbox[t]{1.6in}{Mortality rate of zooplankton} & ${M_{z}}$   &    ${0.0145}$  &
$day^{-1}$ &\parbox[t]{1.55in}{Ray et al., (2001)} \\[2.5ex]
 \hline
\parbox[t]{1.6in}{Excretion of fish} & ${ E_{f} }$  &\parbox[t]{0.3in}{ $0.049$}&
$day^{-1}$
 &\parbox[t]{1.55in}{ Mandal et al., (2012)} \\[2.5ex]
 \hline
 \parbox[t]{1.6in}{Mortality  rate of fish} & ${ M_{f} }$  &\parbox[t]{0.3in}{ $0.021$}&
$day^{-1}$
 &\parbox[t]{1.55in}{Mandal et al., (2012)} \\[2.5ex]
 \hline
 \parbox[t]{1.6in}{Respiration rate of Fish} & ${ r_{f} }$  &\parbox[t]{0.3in}{ $0.0125$}&
$day^{-1}$
 &\parbox[t]{1.55in}{ Mandal et al., (2012)} \\[2.5ex]
 \hline
 \parbox[t]{1.6in}{Harvest rate of Carnivorous Fish} & ${ H_{f} }$  &\parbox[t]{0.3in}{ $0.1090$}&
$day^{-1}$
 &\parbox[t]{1.55in}{ Mandal et al., (2012)} \\[2.5ex]
 \hline
 \parbox[t]{1.6in}{Growth rate of Carnivorous Fish} & ${ g_{f} }$  &\parbox[t]{0.3in}{ $0.6894$}&
$day^{-1}$
 &\parbox[t]{1.55in}{ Mandal et al., (2012)} \\[2.5ex]
 \hline
 \parbox[t]{1.6in}{Half saturation constant for zooplankton grazing by Carnivorous Fish} & ${ k_{F} }$  &\parbox[t]{0.3in}{ $10$}&
\parbox[t]{1.1in}{dimensionless}
 &\parbox[t]{1.55in}{ Mandal et al., (2012)} \\[2.5ex]
 \hline
 \parbox[t]{1.6in}{Half saturation constant for phytoplankton grazing by zooplankton} & ${ k_{Z} }$  &\parbox[t]{0.3in}{ $38$}&
\parbox[t]{1.1in}{dimensionless}
 &\parbox[t]{1.55in}{ Mandal et al., (2012)} \\[1.5ex]
 \hline
\end{tabular}\label{Table:t_1}\\[3.5ex]
\end{table}

\begin{figure}
\centering
\begin{tabular}{c}\label{f1}
(a) \epsfig{file=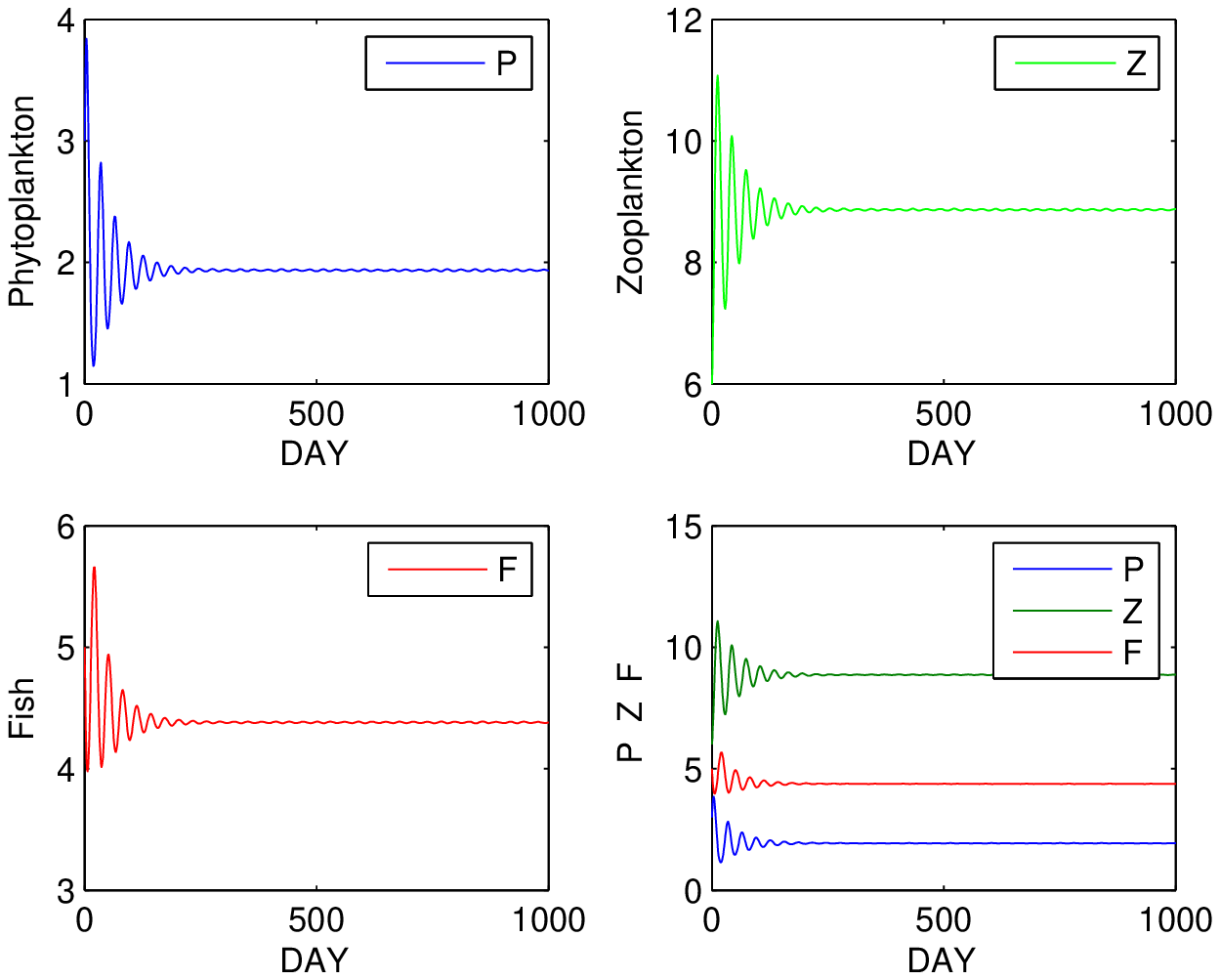,width=0.9\linewidth,height=6cm,clip=} \\
(b) \epsfig{file=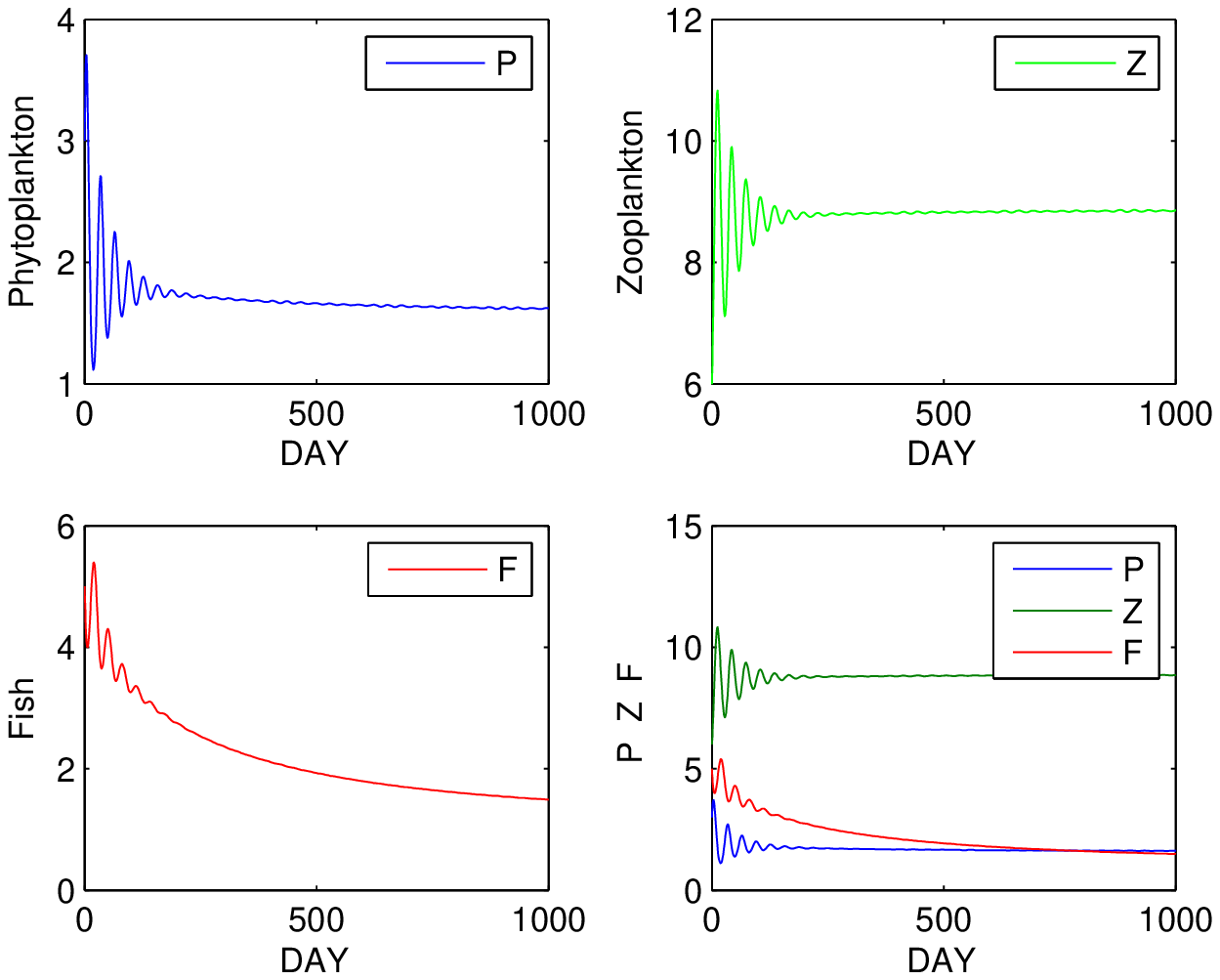,width=0.9\linewidth,height=6cm,clip=}\\
(c)\epsfig{file=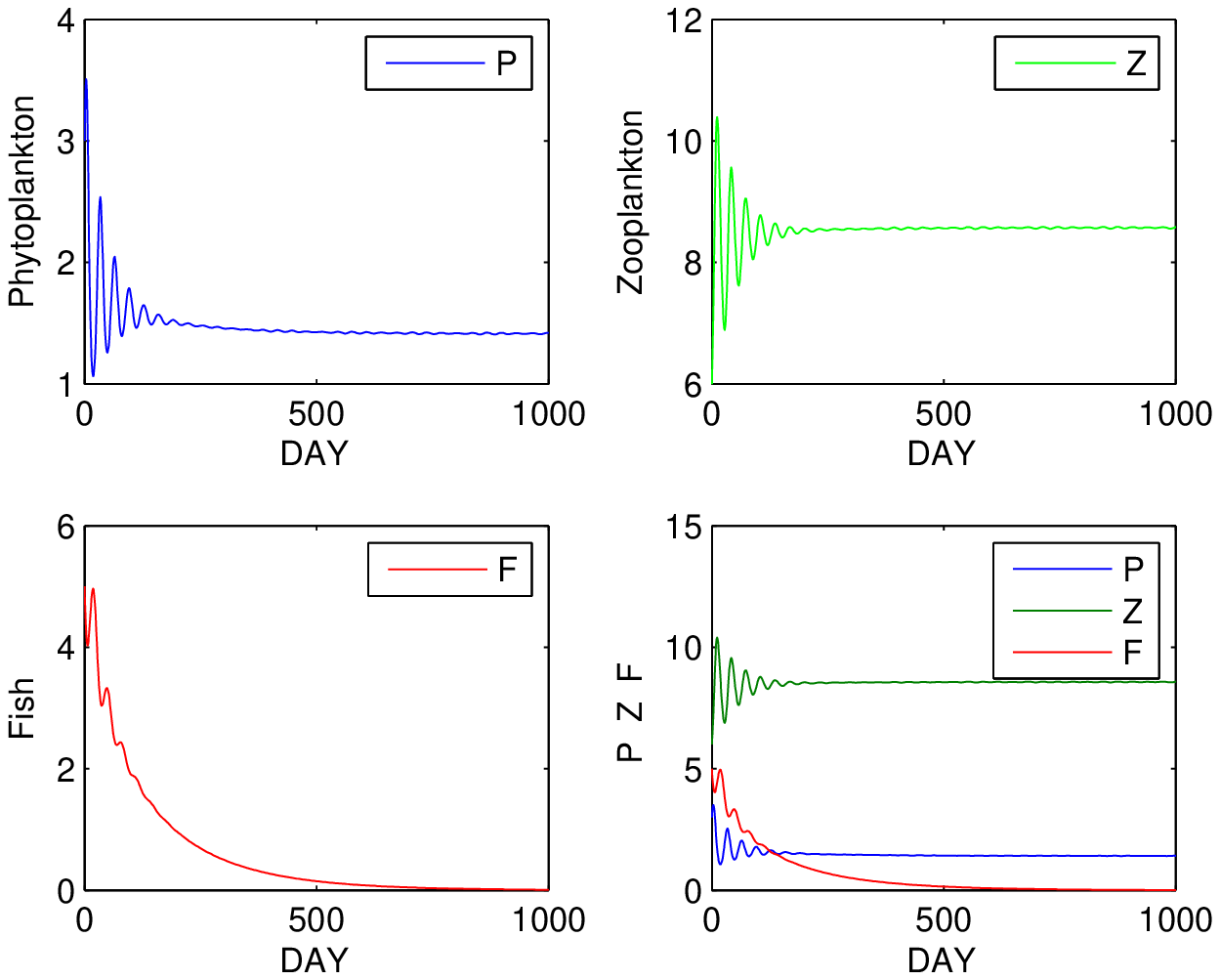,width=0.9\linewidth,height=6cm,clip=}\end{tabular}
\caption{\textrm{\small  Global stability of the PZF dynamical system for different parameter conditions. The graphics depicted in this Figure are obtained by varying the downstream and upstream salinity effect. (a) salinity of downstream $(s_{d}) =12.30,$ salinity of upstream $(s_{u}) = 8.23$; (b) $s_{d} =12.30,$ $s_{u} = 8.33;$ (c) $s_{d} =12.30, s_{u} = 8.51.$}}\end{figure}

\begin{figure}
\begin{center}
(a)\includegraphics[width=18cm, height = 6cm]{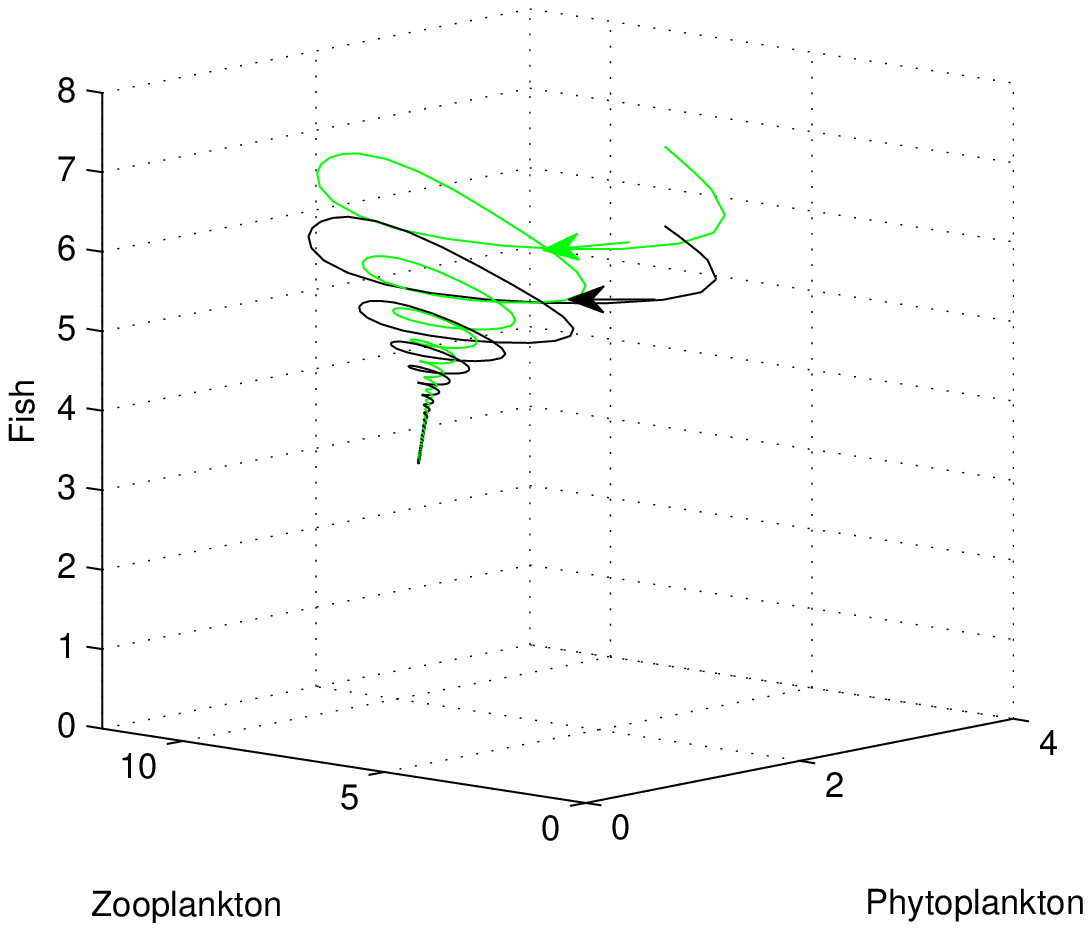}
(b)\includegraphics[width=18cm, height = 6cm]{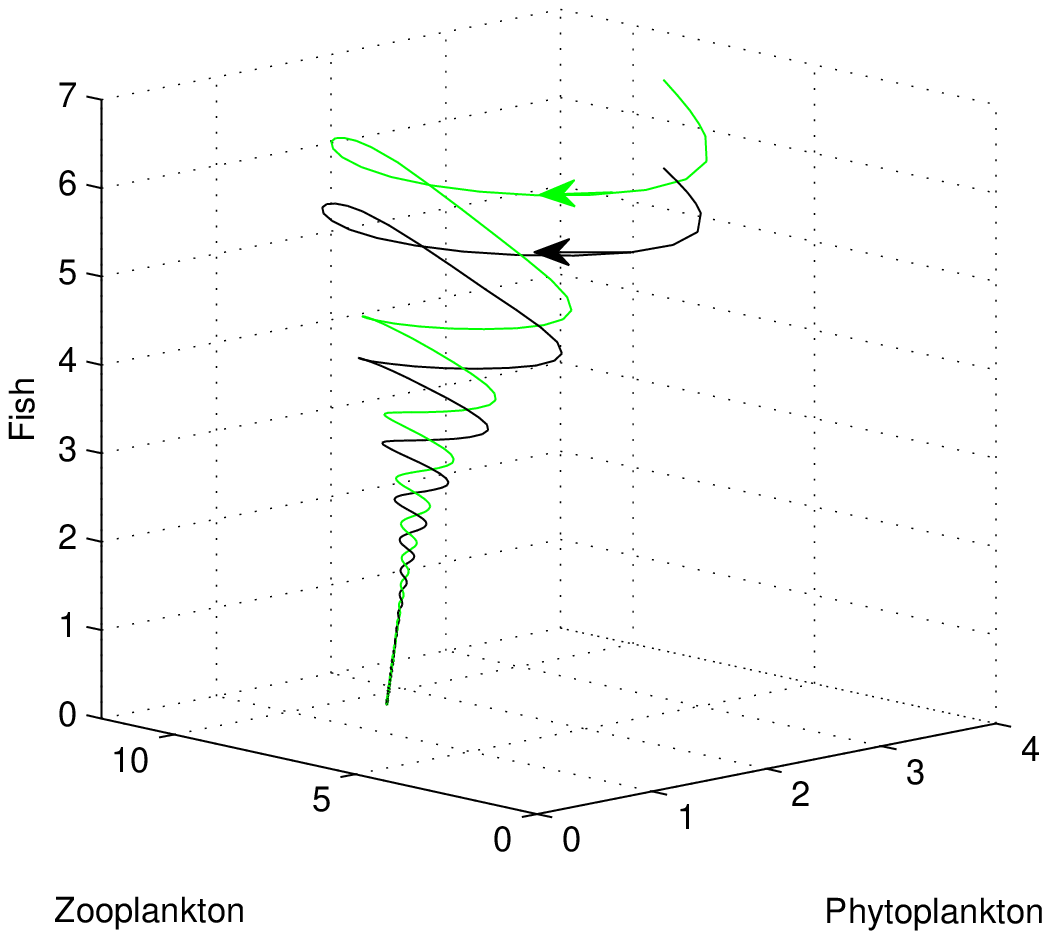}
(c)\includegraphics[width=18cm, height = 6cm]{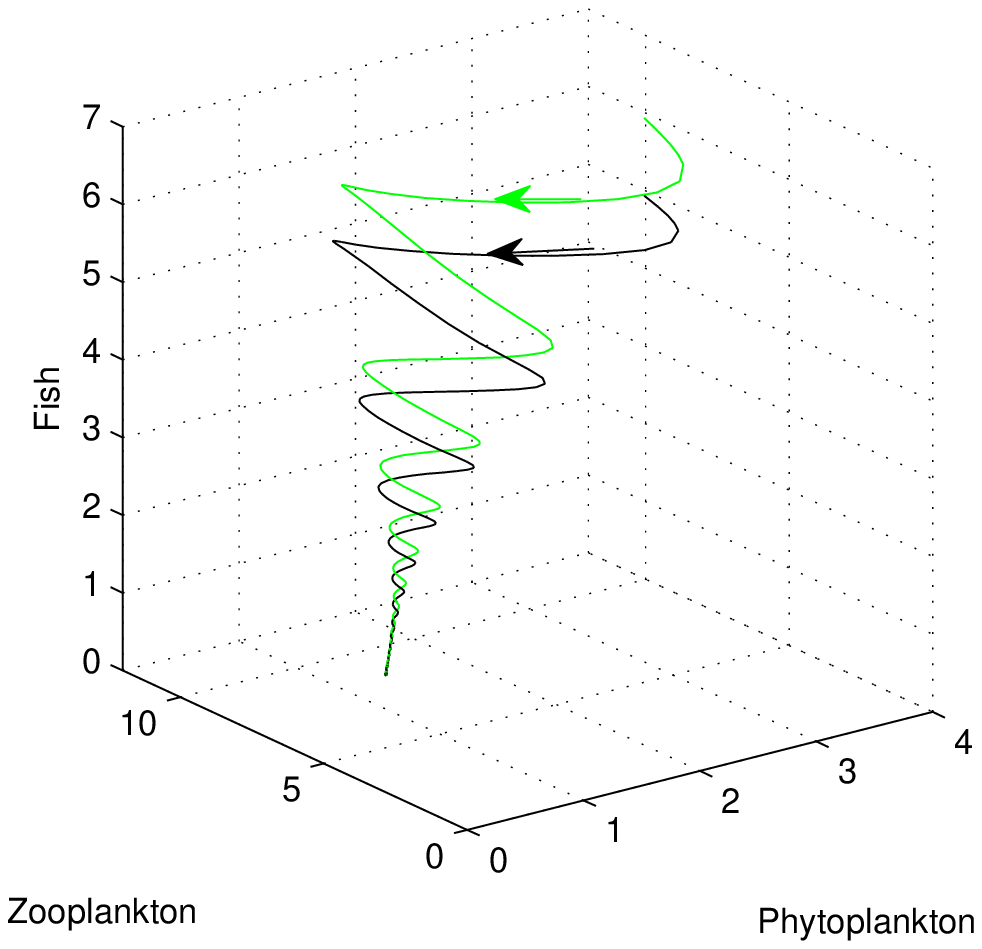}
\end{center}
\caption{\textrm{\small  3D view of global stability of the PZF system for different parameter conditions with the same parameter values used in Figure 1. (a) salinity of downstream $(s_{d}) =12.30,$ salinity of upstream $(s_{u}) = 8.23$; (b)$s_{d} =12.17,$ $s_{u} = 8.33;$ (c) $s_{d} =12.30, s_{u} = 8.51.$}}
\end{figure}

\begin{figure}
\begin{center}
(a)\includegraphics[width=12cm, height = 9cm]{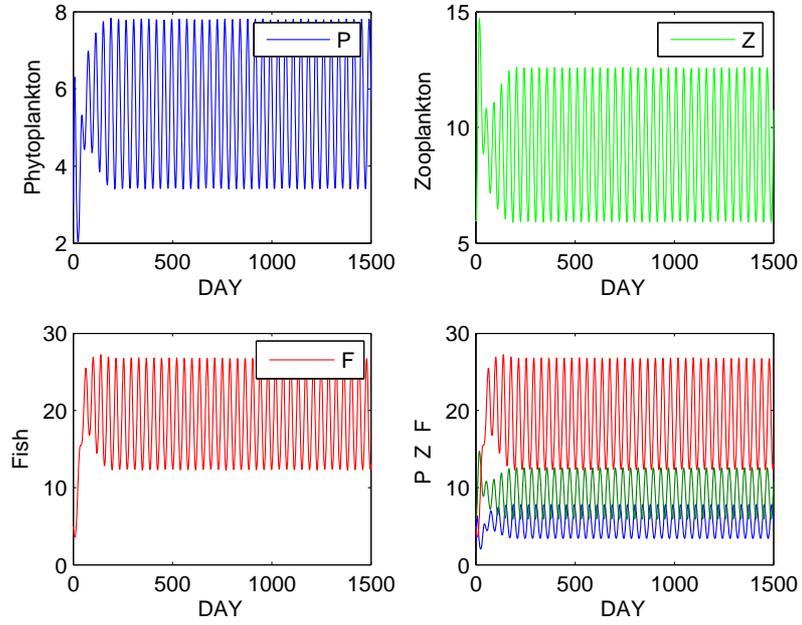}
(b)\includegraphics[width=12cm, height = 9cm]{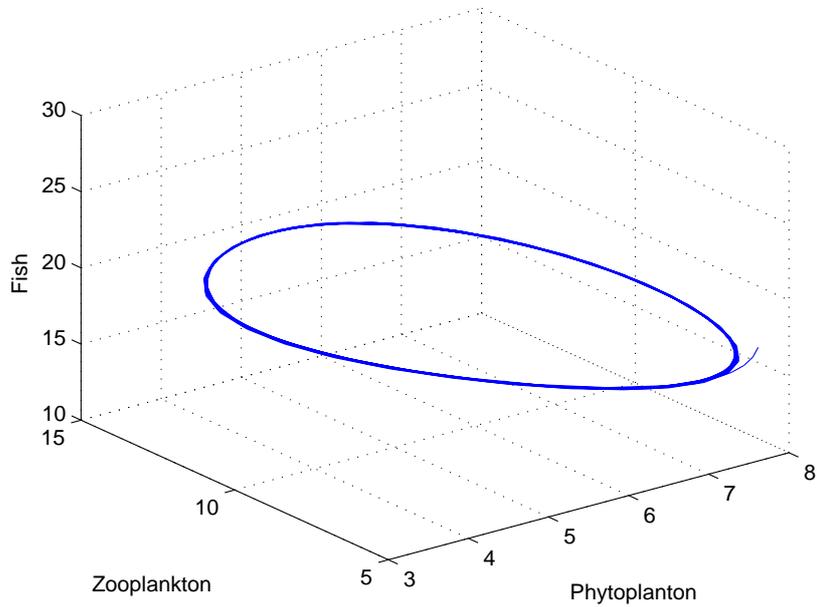}
\end{center}
\caption{\textrm{\small  Limit cycle behaviour of the PZF dynamics for different parameter conditions. (a) 2D view of limit cycle around $E_{*}$ of the system with salinity of downstream $s_{d} =12.30,$ salinity of upstream $s_{u} = 6.25$; (b) 3D view of limit cycle with $s_{d} =12.17, s_{u} = 6.95,$ and other parameter are same as in Figure 1.}}
\end{figure}

\begin{figure}
\begin{center}
(a)\includegraphics[width=12cm, height = 9cm]{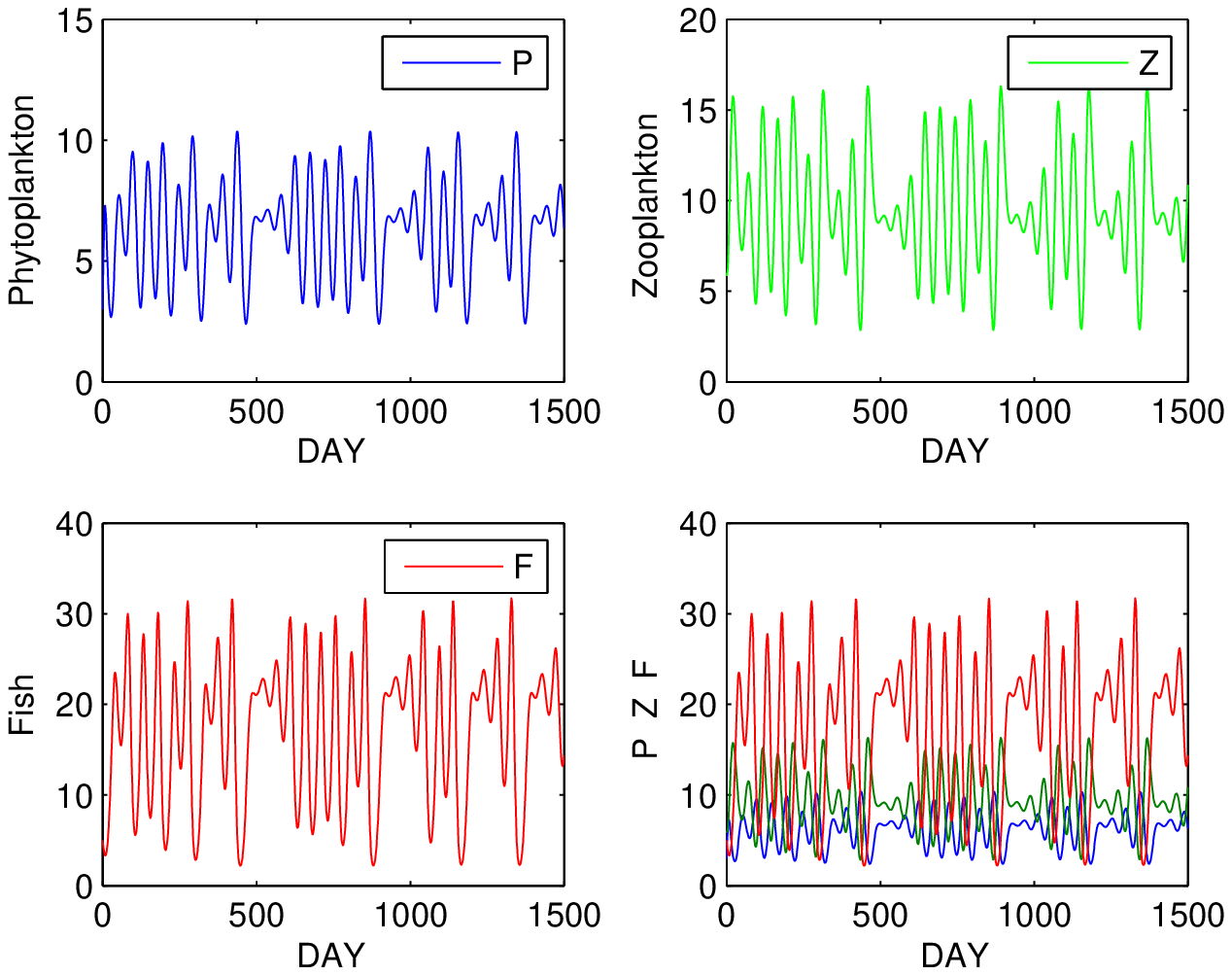}
(b)\includegraphics[width=12cm, height = 9cm]{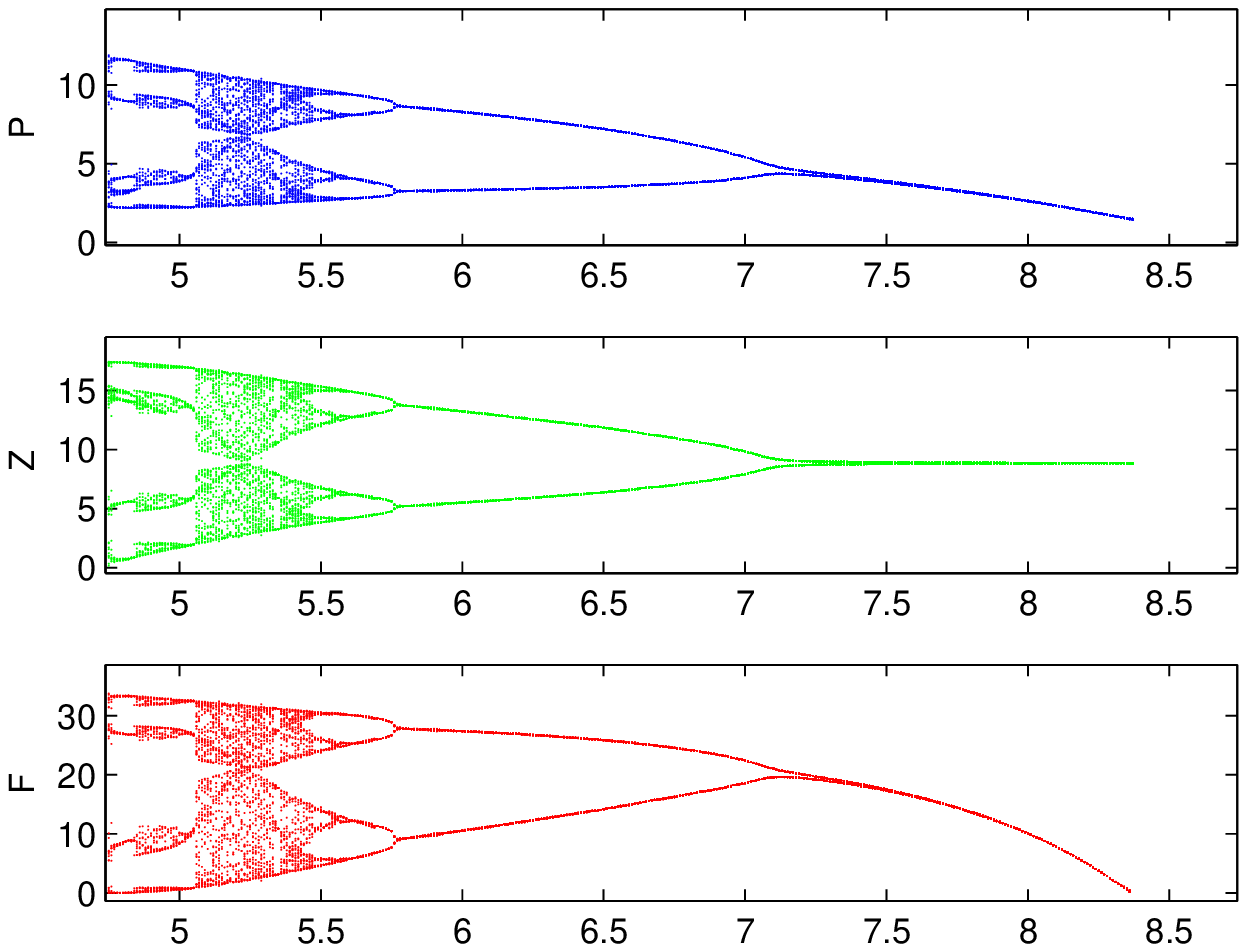}
\end{center}
\caption{\textrm{\small (a) Top panel: represents the chaotic nature of the species with downstream  salinity $s_{d} =12.30,$ and upstream salinity $s_{u} = 5.25$; (b) Bottom panel: represents the same chaotic nature of the species with respect to the variation of $s_{u}$ within the range 5 and 8.5, while $s_{d}$ is fixed at 12.30.}}
\end{figure}

\begin{figure}
\begin{center}
\includegraphics[width=12cm, height = 12cm]{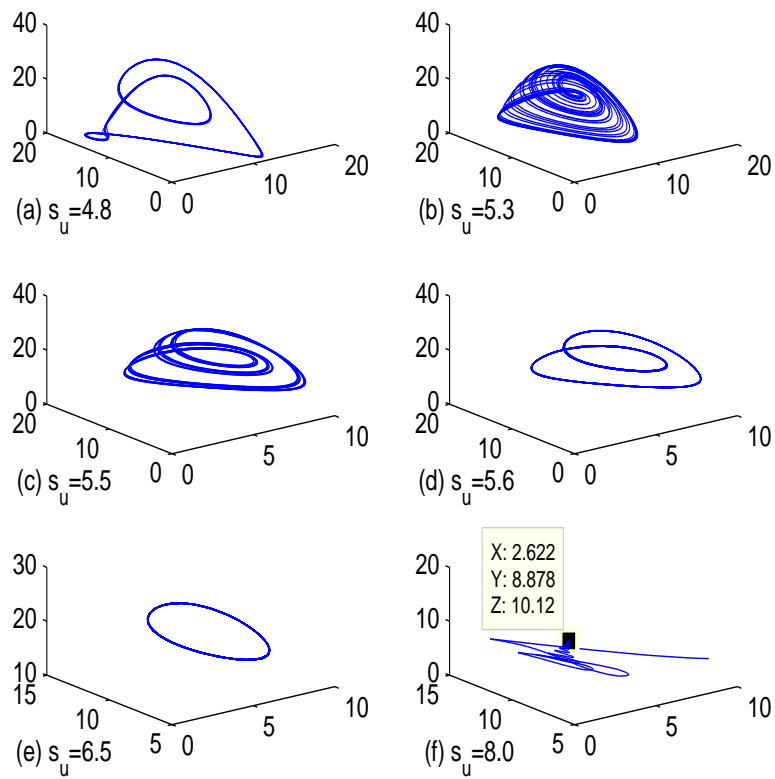}
\end{center}
\caption{\textrm{\small The PZF system executed different phase diagram starting with a periodic orbit of period 6 to a stable solution plot around an equilibrium point $(2.622, 8.878, 10.12)$ via a chaotic orbit, a orbit of period 8, a period doubling and a limit cycle respectively when the salinity of upstream  varies from$s_{u} = 4.50\,~ \text{to}~\, 8.0$ while downstream salinity $s_{d}$ is fixed at 12.30.}}
\end{figure}
\section{Discussion:}
Hooghly-Matla estuarine ecosystem is highly dynamic as the environment changes (change in water temperature and salinity) and as more new species invade (cf. Rosith, (2013) and Bhaumik, (2013)). This estuarine system is mixohaline in nature, where the salinity ranges from freshwater condition ($< 2$ppt) to 22 ppt at various points along the stretch of the estuary (cf. Mandal et al., (2012)).  In the pre-monsoon, fresh water runoff from the upstream becomes very less and due to tidal influence of the adjacent Bay of Bengal, the salinity increases. Throughout the year, a gradient of salinity is observed between upstream and downstream area of the estuary (cf. Mandal et al., (2009)). Previous studies indicate that changes in the estuarine fish assemblage are regulated by associated changes in the salinity and estuarine mouth morphology (cf. Gillanders et al., (2011)). The Hooghly-Matla estuarine system has funnel shaped sea face which is appropriate for optimum tidal flux (cf. Rosith, (2013) and Bhaumik, (2013)).\\
Several research works are published on fish depletion in the Hooghly-Matla estuary. It is observed that there is gradual regime shift of saline water from the upper stretches of estuary to the downstream of the estuary during the past three decades and it happened when of Farakka Barrage is commissioned in 1975 (cf. Sinha, (1996)). The shift of estuarine zone towards the downstream is due to the increased ingress of fresh water to the Hooghly-Matla estuary, which resulted in the extension of freshwater zone. This leads to disappearance of stenohaline fishes from the system. Barron et al., (2002) reported that phytoplankton species of estuarine origin are more tolerant to low salinity than oceanic species. This study suggests the co-existence of euryhaline fishes and estuarine phytoplankton in the Hooghly-Matla estuarine system. The present study is in agreement with the above findings.\\
Another study which supports the present research is the ionic balance of zooplankton in estuarine environment. It is observed that the propensity of invertebrate species richness near the upstream of the estuary decreases as the salinity reaches the critical values (5 to 8 ppt). It occurs due to the inability of invertebrates to regulate specific ionic concentrations at and below the critical salinity (cf. Khlebovich, (1969)). Freshwater ingress in the monsoon period changes the salinity and this in turn indirectly effects the dynamics of fish population of the estuary.\\
The present situation of Hooghly-Matla estuarine system represents stressed condition of the fish species in the estuary. Model shows that the existence of fish in the system is possible only when the growth rate of carnivorous fish population is greater than cumulative effect of excretion rate, respiration rate ($g_f > m_3$). The growth rate is dependent on predation of zooplankton by carnivorous fish population. Energy allocation is dependent on the life-history strategy of a fish species. It is widely accepted that both material and energy are mobilized and reallocated for reproduction in fishes (cf. Jobling, (1995)). This indicates that there exists energy trade off in physiological processes of the estuarine species, where most of the energy earned through predation is allocated for reproduction rather respiration and excretion. Hence, there is growth of fish population in the estuary. We have also shown that the salinity effect has a prominent role in significantly stabilizing the coexistence equilibrium in our model (cf. Figure 1 and Figure 2).

{\bf Acknowledgement:} Mr. Ujjwal Roy and Dr. N. C. Majee are thankful to the Department of Mathematics, Jadavpur University for providing all the facilities to complete this research work. The corresponding author  Dr. S. Sarwardi is thankful to the Department of Mathematics, Aliah University for providing opportunities to perform the present work. He is also thankful to his Ph.D. supervisor Prof. Prashanta Kumar Mandal, Department of Mathematics, Visva-Bharati for his generous help and continuous encouragement while preparing this manuscript. Prof. Santanu Ray is thankful to the Department of Zoology, Visva-Bharati University (a Central University) for opportunities to perform the present work.

\end{document}